\documentclass{amsproc}

\usepackage{amsmath, amsfonts, amsthm, mathrsfs, amssymb}

\usepackage{graphicx}  
\usepackage{tikz}
\usetikzlibrary{cd, trees}
\usepackage{tikz-qtree}
\usepackage{forest}

\usepackage{algorithm}
\usepackage{algorithmic}
\usepackage[backend=biber]{biblatex}
\usepackage{url}

\newtheorem{theorem}{Theorem}[section]
\newtheorem{lemma}[theorem]{Lemma}
\newtheorem{proposition}[theorem]{Proposition}
\newtheorem{corollary}[theorem]{Corollary}

\theoremstyle{definition}
\newtheorem{definition}[theorem]{Definition}
\newtheorem{example}[theorem]{Example}

\theoremstyle{remark}
\newtheorem{remark}[theorem]{Remark}

\numberwithin{equation}{section}
\begin{document}

\title{Algorithm for Motivic Hilbert Zeta Function of Some Curve Singularities}




\author{Yizi Chen}
\address{ETH Zurich, Zurich, 8093, Switzerland}
\email{yizi.chen@ethz.ch}
\thanks{}

\author{Hussein Mourtada}
\address{Universit\'e Paris Cit\'e, Sorbonne Universit\'e \\CNRS, IMJ-PRG, F-75013 Paris}
\email{hussein.mourtada@imj-prg.fr}
\thanks{}


\author{Wenhao Zhu}
\address{Universit\'e Paris Cit\'e,  IMJ-PRG, F-75013 Paris}
\email{wenhao.zhu.alg@gmail.com}
\date{}

\begin{abstract}
We develop algorithms to compute two versions of the motivic Hilbert zeta function for curve singularities: the classical version, applicable to singularities with a monomial valuation semigroup or to singular curves defined by \(y^{k}=x^{n}\) with \(\gcd(k,n)=1\), and a finer version introduced by the first and third authors together with Mounir Hajli, which currently applies to the specific family \(y^{k}=x^{n}\) where \(\gcd(k,n)=1\).

It is well known that the Hilbert scheme of points on a smooth curve is isomorphic to the symmetric product of the curve. However, the geometry of the Hilbert scheme of points on singular curves remains much less understood. Our algorithms compute the motivic Hilbert zeta functions
\[
Z_{(C,O)}^{\mathrm{Hilb}}(q) \in K_{0}(\mathrm{Var}_{\mathbb{C}})[[q]], \qquad 
Zm_{(C,O)}^{\mathrm{Hilb}}(a^{2},q^{2}) \in K_{0}(\mathrm{Var}_{\mathbb{C}})[[a^{2}, q^{2}]],
\]
for such curve singularities, expressed as formal power series with coefficients in the Grothendieck ring of complex varieties.

The main computational difficulty arises from the fact that \(\Gamma\) is infinite. To overcome this, we approximate \(\Gamma\) by truncating it to a suitable finite subset, which allows the algorithms to run effectively. We analyze the time complexity of the method and provide an estimate for the effective finite length of \(\Gamma\) required to obtain reliable results. A Python implementation of the algorithms is available at \url{https://github.com/whaozhu/motivic_hilbert.}
\end{abstract}

\maketitle

\section{Introduction}
Given a complex curve singularity $(C,O)$ and an integer $\ell \in \mathbb{Z}_{\geq 0}$, the $\ell$-th punctual Hilbert scheme of $(C,O)$, denoted $C^{[\ell]}$, is the moduli space parametrizing $0$-dimensional subschemes of $(C,O)$ of length $\ell$. This is a special case of Grothendieck's Hilbert scheme, which more generally parametrizes subschemes of projective space with fixed Hilbert polynomial.\vskip 0.1cm
From an algebraic viewpoint, let $A = \mathcal{O}_{C,O}$ denote the local ring of $C$ at $O$. Then the $\mathbb{C}$-points of $C^{[\ell]}$ are in bijection with ideals of colength $\ell$ in $A$:
\[
C^{[\ell]} := \left\{ I \subset A \mid I \text{ is an ideal and } \dim_\mathbb{C} \frac{A}{I} = \ell \right\}.
\]\vskip 0.1cm

Let $m(I)$ be the minimal number of generators of ideal $I$. For $m\in \mathbb{Z}_{\geq 1}$, denote the subvariety of punctual Hilbert scheme parametrizing ideals of $A$ with $m$ minimal number of generators:  
\[C^{[\ell],m}:=\{I\in C^{[\ell]}\mid m(I)=m\}.\]

The goal of this paper is to develop  algorithms to compute the motivic Hilbert zeta function 
\[Z_{(C, O)}^{\text {Hilb}}(q) :=  \sum_{\ell=0}^{\infty} [C^{[\ell]}] q^{\ell} \in K_{0}(Var_{\mathbb{C}})[[q]],\]
for a germ of a singular curve $(C,O)$ with a monomial valuation semigroup in the sense of G. Pfister and J.H.M. Steenbrink \cite{pfister1992reduced}, as well as for singular curves defined by $y^{k}=x^{n}$, and to compute the new motivic Hilbert zeta function introduced by the first and third authors together with Mounir Hajli \cite{hajli2025geometry}, 
\[
Zm_{(C,O)}^{\mathrm{Hilb}}(a^2,q^2)= \sum_{\ell \geq 0} \sum_{m \geq 1} q^{2\ell} (1 - a^2)^{m-1} [C^{[\ell],m}]  \in K_{0}(\mathrm{Var}_{\mathbb{C}})[[a^2, q^2]],
\] 
for singular curves defined by $y^{k}=x^{n}$.

In Section \ref{section: Punctual Hilbert scheme and monomial semigroup}, we recall results on the punctual Hilbert scheme due to G.~Pfister and J.~H.~M.~Steenbrink \cite{pfister1992reduced}.  
Let \(\overline{A}\) denote the integral closure of \(A\) and consider the delta invariant \(\delta := \dim_k \overline{A}/A\). Note that \(\overline{A}\) is a discrete valuation ring; the induced valuation is denoted by \(v \colon \overline{A} \longrightarrow \mathbb{N}\). Let \(\Gamma := v(A \setminus \{0\})\) be its semigroup. It is shown in \cite{pfister1992reduced, piontkowski2007topology} that the Hilbert scheme \(C^{[\ell]}\) can be embedded into a Grassmannian variety \(\operatorname{Gr}(\delta,2\delta)\), where \(\delta\) is the \(\delta\)-invariant of \((C,O)\).  From this embedding one obtains equations that define \(C^{[\ell]}\) (Theorem 3 in \cite{pfister1992reduced}). Nevertheless, studying the geometry of \(C^{[\ell]}\) directly via these equations is difficult because one must consider the image of \(C^{[\ell]}\) under the Plücker embedding  
\(
\operatorname{Gr}(\delta,2\delta) \hookrightarrow \mathbb{P}\!\Bigl(\bigwedge\nolimits^{\delta} \mathbb{C}^{2\delta}\Bigr),
\)  
and such computations become very involved.

In \cite{oblomkov2018hilbert},  A. Oblomkov, J. Rasmussen, and V. Shende proved that for each \(\ell \geq 1\) the scheme \(C^{[\ell]}\) admits a stratification (see  \eqref{stratification of Hilbert scheme}) and they give defining equations for each stratum (Proposition 12 in \cite{oblomkov2018hilbert}). Given a subsemimodule \(\Delta\) of \(\Gamma\) satisfying \(\Gamma + \Delta \subseteq \Delta\), the corresponding stratum parametrizes ideals whose valuation set equals \(\Delta\); it is denoted by \(C^{[\Delta]}\):
\[
C^{[\Delta]} := \left\{ I \in C^{[\ell]} \mid \ell = \#(\Gamma \setminus \Delta),\ v(I \setminus \{0\}) = \Delta \right\}.
\]

In Section \ref{section: Tree structure associated with valuation semigroup}, we recall the work in \cite{hajli2025geometry}. Inspired by the recent work of Y. Soma and M. Watari \cite{soma2014punctual} on the geometry of Hilbert schemes, the authors associate to a curve singularity with semigroup $\Gamma$ a \emph{leveled graph} $G_\Gamma$, which is defined as follows. The vertices at level $\ell$ are the elements of the set
\[
\mathscr{D}_\ell := \left\{ \Delta \subset \Gamma \mid \Delta \text{ is a $\Gamma$-subsemimodule with } \#(\Gamma \setminus \Delta) = \ell \right\}.
\]
There is a directed edge from $\Delta \in \mathscr{D}_\ell$ to $\Delta' \in \mathscr{D}_{\ell-1}$ if 
\[
\Delta' = m(\Delta) := \Delta \cup \{\gamma_\Delta\}, \quad \text{where} \quad \gamma_\Delta := \max(\Gamma \setminus \Delta).
\]
The authors study the graph $G_\Gamma$ and prove that it has the structure of a tree.

In Section \ref{section: Piecewise fibrations}, we recall results from \cite{hajli2025geometry}, which assert that when $C$ is a curve with a monomial valuation semigroup or is defined by $y^{k}=x^{n}$, an edge of the tree $G_{\Gamma}$ connecting two vertices $\Delta$ and $m(\Delta)$ corresponds to a geometric morphism,  a piecewise fibration over its image,  between $C^{[\Delta]}$ and $C^{[m(\Delta)]}$ (Theorem~\ref{piecewise fibration}). This morphism encodes the relationship between the classes of these strata in the Grothendieck ring (Corollary~\ref{mononial case,relation in grothendieck}).

In Section~\ref{Algorithm 1}, we present Algorithm~\ref{algorithm motivic} for computing the motivic Hilbert zeta function of monomial curve singularities or curve defined by $y^p = x^q$. The computation relies on several auxiliary algorithms (Algorithms~\ref{algorithm min generators}--\ref{algorithm tree}), each performing a specific task: identifying minimal generators, computing syzygies, and generating sets of $\Gamma$-subsemimodules from the valuation semigroup. To ensure computational feasibility, all routines work with a truncated set $\Gamma_{\leq n}$. The time complexity of the complete procedure is analyzed in detail.

In  Section~\ref{section Algorithm for Zm}, we  introduce Algorithm~\ref{new algorithm motivic}, Algorithm~\ref{algm: C^Delta,<m}, and Algorithm~\ref{algm: C^Delta,i} to compute the \emph{new} motivic Hilbert zeta function $Zm_{(C, O)}^{\text {Hilb}}$ for curve singularities $C$ defined by $y^{p}=x^{q}$, where $\gcd(p,q)=1$.

In Section~\ref{Effective finite length}, we introduce $N$ as the minimal effective truncation bound for $\Gamma$ required by Algorithm~\ref{algorithm motivic} and provide an estimate for its range (Theorem~\ref{limite of effective bound}).

In Section~\ref{section Experimentation}, we focus on curve singularities of type $W_{8}$, i.e., those whose complete local rings are isomorphic to $\mathbb{C}[[t^{4}, t^{5}, t^{6}]]$. The geometry of punctual Hilbert schemes for certain ADE curve singularities (types $A_{2k}$ with $k\geq 1$, $E_6$, and $E_8$) has been studied by Y.~Soma \cite{soma2015punctual} and M.~Watari \cite{soma2014punctual}. More recently, Watari computed the motivic Hilbert zeta function for these types \cite{watari2024motivic}. Our Algorithm~\ref{algorithm motivic} provides an alternative way to compute the motivic Hilbert zeta function for such specific curve singularities (see examples at \url{https://github.com/whaozhu/motivic_hilbert}).

\section{Punctual Hilbert schemes and monomial semigroup} \label{section: Punctual Hilbert scheme and monomial semigroup}

In this section, we recall the definition of the punctual Hilbert schemes of a curve singularity and review some fundamental properties. Let $(C,O)$ be the germ of a unibranch curve singularity defined over $\mathbb{C}$, and let $A := \mathcal{O}_{C,O}$ denote its local ring at $O$. Then $A$ is a complete one-dimensional local domain, and its integral closure $\overline{A}$ in its field of fractions is a discrete valuation ring, isomorphic to $\mathbb{C}[[t]]$.

Let $v \colon \overline{A} \setminus \{0\} \to \mathbb{N}$ be the associated discrete valuation, extended by $v(0) = \infty$, which maps a nonzero power series in $\mathbb{C}[[t]]$ to its order in $t$. The \emph{value semigroup} of $A$ (or of the singularity $(C,O)$) is defined as
\[
\Gamma := v(A \setminus \{0\}) \subseteq \mathbb{N}.
\]
This is a numerical semigroup: a submonoid of $\mathbb{N}$ with finite complement.

For $n \in \mathbb{N}$, define the ideal
\[
\overline{I}(n) := \{ f \in \overline{A} \mid v(f) \geq n \}
\]
in $\overline{A}$, and set $I(n) := \overline{I}(n) \cap A$, which is an ideal in $A$. The \emph{conductor} of $A$ in $\overline{A}$ is the smallest integer $c$ such that $\overline{I}(c) \subset A$, i.e.,
\[
c := \min\left\{ n \in \mathbb{N} \mid \overline{I}(n) \subset A \right\}.
\]
Equivalently, $c$ is the smallest element such that $t^c \mathbb{C}[[t]] \subseteq A$ under a fixed isomorphism $\overline{A} \cong \mathbb{C}[[t]]$.

The \emph{$\delta$-invariant} of the singularity is defined as
\[
\delta := \dim_{\mathbb{C}}(\overline{A}/A) = \#(\mathbb{N} \setminus \Gamma).
\]
It measures the arithmetic genus defect of the singularity. We have the inequalities
\[
\delta + 1 \leq c \leq 2\delta,
\]
and $c = 2\delta$ if and only if $A$ is Gorenstein~\cite{serre1959groupes}.

For $\ell \in \mathbb{Z}_{>0}$, the \emph{$\ell$-th punctual Hilbert scheme} of $(C,O)$ is the moduli space
\[
C^{[\ell]} := \left\{ I \subset A \mid I \text{ is an ideal and } \dim_\mathbb{C} \frac{A}{I} = \ell \right\}.
\]
It parametrizes zero-dimensional subschemes of length $\ell$ supported at $O$.\\

Pfister and Steenbrink~\cite{pfister1992reduced} showed that the punctual Hilbert schemes of a unibranch curve singularity can be embedded as closed subvarieties of a Grassmannian. More precisely, let $\mathscr{M}$ denote the reduced subscheme of the Grassmannian $\operatorname{Gr}(\delta, \overline{A}/I(2\delta))$ defined by
\[
\mathscr{M} := \left\{ W \in \operatorname{Gr}(\delta, \overline{A}/I(2\delta)) \mid W \text{ is an } A\text{-submodule of } \overline{A}/I(2\delta) \right\}.
\]
One can verify that $\mathscr{M}$ is a linear subvariety of the Grassmannian—i.e., it is defined by a system of linear equations in the Plücker coordinates.

\begin{proposition}[\cite{pfister1992reduced}, Theorem~3]\label{prop:pfister-steenbrink}
For every $\ell > 0$, there exists a closed embedding
\[
\phi_\ell \colon C^{[\ell]} \longrightarrow \mathscr{M}.
\]
Moreover, $\phi_\ell$ is bijective (as a map of sets) when $\ell \geq c$, the conductor of the singularity.
\end{proposition}

Since $c \leq 2\delta$, this implies that for $\ell \geq 2\delta$, the structure of $C^{[\ell]}$ stabilizes in the sense that its points are completely determined by the fixed ambient variety $\mathscr{M}$. Consequently, to understand the classes $[C^{[\ell]}]$ in the Grothendieck ring of varieties $K_0(\mathrm{Var}_k)$, it suffices to study the punctual Hilbert schemes for $\ell \leq 2\delta$.

A powerful method to compute these classes is to stratify each $C^{[\ell]}$ into constructible subsets using the valuation $v$ and the value semigroup $\Gamma = v(A \setminus \{0\})$.   A subset $\Delta$ of $\Gamma$ is called a \emph{$\Gamma$-subsemimodule} if $\Delta + \Gamma \subset \Delta$. For an ideal $I$, we denote $\Gamma(I):= v(I\setminus \{0\})$, which is a $\Gamma$-subsemimodule. 

For any $\Gamma$-subsemimodule $\Delta$, we define the stratum
\[
C^{[\Delta]} := \left\{ I \in C^{[\ell]} \mid \ell = \#(\Gamma \setminus \Delta),\ v(I \setminus \{0\}) = \Delta \right\}.
\]

The following lemma, taken from \cite{pol2016singularites}, will be essential.

\begin{lemma}[\cite{pol2016singularites}]\label{ideal codimension}
For a positive integer $\ell$, an ideal $I$ of $A$ belongs to $C^{[\ell]}$ if and only if $\#\bigl(\Gamma \setminus \Gamma(I)\bigr)=\ell$.
\end{lemma}

Define \[\mathscr{D}_{\ell}:=\{\Delta\subset \Gamma \mid \Delta \text{ is a $\Gamma$-subsemimodule with } \#(\Gamma\setminus \Delta)=\ell \},\] the set of $\Gamma$-subsemimodules whose complement in $\Gamma$ has cardinality $\ell$. Lemma~\ref{ideal codimension} then yields the stratification
\begin{equation}\label{stratification of Hilbert scheme}
C^{[\ell]} = \bigsqcup_{\Delta \in \mathscr{D}_{\ell}} C^{[\Delta]}.
\end{equation}

\begin{remark}
With the notation above, $C^{[\Delta]}$ coincides with the intersection of $C^{[\ell]}$ with the Schubert cell of the Grassmannian $\operatorname{Gr}(\delta, A/I(2\delta))$ (see Lemma~5 in \cite{pfister1992reduced}). Hence $C^{[\Delta]}$ is a locally closed subset of $C^{[\ell]}$.
\end{remark}

In the following section, we will describe a method for determining all possible elements of $\mathscr{D}_{\ell}$.

To end this section, we introduce the notion of monomial semigroups in the sense of G.~Pfister and J.~H.~M.~Steenbrink \cite{pfister1992reduced}.

\begin{definition}[\cite{pfister1992reduced}]
A \emph{monomial curve singularity} over a field $k$ is an irreducible curve singularity whose local ring is isomorphic to $A = k[[t^{a_{1}},\dots,t^{a_{m}}]]$ with $\gcd(a_{1},\dots,a_{m})=1$.
\end{definition}

In this case we may assume that the semigroup $\Gamma$ is generated by its minimal set of generators:
\[
\Gamma = \langle a_{1},\dots,a_{m}\rangle.
\]

\begin{definition}[\cite{pfister1992reduced}]
A semigroup $\Gamma\subset\mathbb{N}$ is called \emph{monomial} if $0\in\Gamma$, $\#(\mathbb{N}\setminus\Gamma)<\infty$, and every reduced irreducible curve singularity with semigroup $\Gamma$ is a monomial curve singularity.
\end{definition}

\begin{proposition}[\cite{pfister1992reduced}]\label{monomial semigroup}
For a semigroup $\Gamma\subset\mathbb{N}$ the following are equivalent:
\begin{enumerate}
\item $\Gamma$ is a monomial semigroup.
\item $\Gamma$ is of one of the following forms:
\begin{enumerate}
\item $\Gamma_{m,s,b}:=\{im \mid i=0,1,\dots,s\}\cup [sm+b,\infty)$ with $1\leq b<m$, $s\geq 1$;
\item $\Gamma_{m,r}:=\{0\}\cup [m, m+r-1]\cup [m+r+1,\infty)$ with $2\leq r\leq m-1$;
\item $\Gamma_{m}:=\{0,m\}\cup [m+2, 2m]\cup [2m+2,\infty)$ with $m\geq 3$.
\end{enumerate}
\item $0\in\Gamma$, $\#(\mathbb{N}\setminus\Gamma)<\infty$, and the following condition holds:\\
if $x\in \mathbb{N}\setminus \Gamma$ and $c(x):=\min\{n\in \mathbb{N}\mid [n,\infty)\subset \Gamma\cup(x+\Gamma)\}$, then
\[
\Gamma\cap (x+\Gamma)\subset [c(x),\infty).
\]
\end{enumerate}
\end{proposition}

\begin{corollary}[\cite{pfister1992reduced, hajli2025geometry}]
Given a $\Gamma$-subsemimodule $\Delta$, let $\gamma_{1}>\gamma_{2}$ be two generators of $\Delta$. Then for any $\sigma\in\operatorname{Syz}\bigl(\langle\gamma_{1},\gamma_{2}\rangle\bigr)$ we have $\sigma\geq c(\Delta)$.
\end{corollary}

\section{The tree structure of subsemimodules}\label{section: Tree structure associated with valuation semigroup}

Let $\Gamma$ be the value semigroup associated with the germ of an irreducible curve singularity $(C,O)$. The goal of this section is to recall that the set of $\Gamma$-subsemimodules can be endowed with a tree structure \cite{hajli2025geometry}.

For $\Delta \in \mathscr{D}_\ell$, let $\gamma_1, \dots, \gamma_n$ be a minimal system of generators of $\Delta$ as a $\Gamma$-subsemimodule. Up to reordering, we may assume $\gamma_1 < \cdots < \gamma_n$. We use the notation
\[
\Delta = (\gamma_1, \dots, \gamma_n)_{\Gamma} := \bigcup_{i=1}^n (\gamma_i + \Gamma),
\]
where the sum denotes the smallest $\Gamma$-subsemimodule containing all $\gamma_i + \Gamma$. By minimality, we have
\[
\Delta \supsetneq \bigcup_{\substack{j=1 \\ j \neq i}}^n (\gamma_j + \Gamma) \quad \text{for all } i \in \{1, \dots, n\}.
\]

The following result appears in \cite{soma2014punctual}:

\begin{lemma}\cite{hajli2025geometry}\label{lem:del}
Let $\Delta = (\gamma_1, \dots, \gamma_n)_{\Gamma}  \in \mathscr{D}_\ell$, where $\{\gamma_1, \dots, \gamma_n\}$ is a minimal system of generators. Then for every $i \in \{1, \dots, n\}$, the set
\[
\Delta \setminus \{\gamma_i\}
\]
belongs to $\mathscr{D}_{\ell+1}$.
\end{lemma}

For $n \in \mathbb{N}$, define
\[
\mathscr{D}_{\ell,n} := \left\{ \Delta \in \mathscr{D}_\ell \mid \Delta \text{ has exactly } n \text{ minimal generators as a } \Gamma\text{-subsemimodule} \right\}.
\]
Note that we have the disjoint decomposition
\[
\mathscr{D}_\ell = \bigsqcup_{n \geq 1} \mathscr{D}_{\ell,n}.
\]

It follows from Lemma~\ref{lem:del} that for each $i \in \{1, \dots, n\}$, there is a canonical \emph{deletion map}
\[
d_{\ell,i} \colon \mathscr{D}_{\ell,n} \to \mathscr{D}_{\ell+1}, \quad \Delta \mapsto \Delta \setminus \{\gamma_i\},
\]
where $\gamma_i$ is one of the minimal generators of $\Delta$.

A natural and important question is whether every element of $\mathscr{D}_{\ell+1}$ arises in this way—i.e., can all $\Gamma$-subsemimodules of codimension $\ell+1$ be obtained by removing a minimal generator from some $\Gamma$-subsemimodule of codimension $\ell$ ? In other words:
\[
\mathscr{D}_{\ell+1} \overset{?}{=} \bigcup_{n \geq 1} \bigcup_{i=1}^n d_{\ell,i}(\mathscr{D}_{\ell,n}).
\]
The answer is affirmative, as will follow from Proposition~\ref{prop:surjective}, which shows that the deletion maps are jointly surjective. 
\begin{definition}\cite{hajli2025geometry}
For a  $\Gamma-$subsemimodule $\Delta$. We define the \emph{Frobenius element} of $\Delta$ by
\[
\gamma_{\Delta} := \max(\Gamma \setminus \Delta).
\]

\end{definition}

\begin{lemma}\label{lem:delta_union}\cite{hajli2025geometry}
Let $\ell \in \mathbb{N}$ and $\Delta \in \mathscr{D}_{\ell+1}$. Then $\Delta \cup \{\gamma_{\Delta}\} \in \mathscr{D}_\ell$.
\end{lemma}

\begin{lemma}\label{lem:gamma_generator}\cite{hajli2025geometry}
Let $\Delta \in \mathscr{D}_{\ell+1}$. Then the Frobenius element $\gamma_\Delta = \max(\Gamma \setminus \Delta)$ is part of every minimal generating set of the $\Gamma$-subsemimodule $\Delta \cup \{\gamma_\Delta\}$.
\end{lemma}

\begin{proposition}\label{prop:surjective}\cite{hajli2025geometry}
The deletion maps cover all of $\mathscr{D}_{\ell+1}$, i.e.,
\[
\mathscr{D}_{\ell+1} = \bigcup_{n \geq 1} \bigcup_{i=1}^n d_{\ell,i}(\mathscr{D}_{\ell,n}).
\]
In other words, every element of $\mathscr{D}_{\ell+1}$ arises as $\Delta' \setminus \{\gamma_i\}$ for some $\Delta' \in \mathscr{D}_\ell$ and some minimal generator $\gamma_i$ of $\Delta'$.
\end{proposition}

It follows from Lemma~\ref{lem:delta_union} that we can define a key map for the sequel: for $\ell \in \mathbb{N}$ with $1 \leq \ell < c$, set
\[
m_{\ell+1} \colon \mathscr{D}_{\ell+1} \to \mathscr{D}_{\ell}, \quad \Delta \mapsto \Delta \cup \{\gamma_{\Delta}\},
\]
where $\gamma_{\Delta} = \max(\Gamma \setminus \Delta)$ is the Frobenius element of $\Delta$. This map, which "adds back" the largest gap in $\Gamma \setminus \Delta$, will play a central role in the recursive study of strata in the punctual Hilbert schemes.

The following definition introduces one of the main objects of this article. Recall that $c$ denotes the conductor of the curve singularity $(C,O)$.
\begin{definition}\label{def:G_Gamma}\cite{hajli2025geometry}
The \emph{$\Gamma$-subsemimodule graph}, denoted $G_\Gamma = (V, E)$, is the directed, levelled graph defined as follows:
\begin{itemize}
    \item \textbf{Vertices:} For each $\ell$ with $1 \leq \ell \leq c$, we denote the set of the vertices at level  $V_\ell:=\mathscr{D}_\ell$, so that
    \[
    V = \bigsqcup_{\ell=1}^{c} V_\ell.
    \]

    \item \textbf{Edges:} For $2 \leq \ell \leq c$, there is a directed edge from a vertex corresponding to $\Delta \in \mathscr{D}_\ell$ to a vertex corresponding to $\Delta' \in \mathscr{D}_{\ell-1}$ if and only if
    \[
    m_\ell(\Delta) = \Delta'.
    \]
    We denote by $E_\ell$ the set of such edges from level $\ell$ to level $\ell-1$, so that
    \[
    E = \bigsqcup_{\ell=2}^{c} E_\ell.
    \]
\end{itemize}
We refer to $G_\Gamma$ as the \emph{value tree} or \emph{semimodule tree} associated with the singularity $(C,O)$.
\end{definition}

\begin{theorem}\label{Tree}\cite{hajli2025geometry}
The $\Gamma$-subsemimodule graph $G_{\Gamma}$ admits a canonical tree structure. The set $\mathscr{D}_1$ consists of a single element, $\Gamma \setminus \{0\}$, which we designate as the root of the graph $G_{\Gamma}$.
\end{theorem}

\begin{example}\label{ex:34}
Let $A = k[[t^3, t^4]]$, the local ring of the plane curve singularity defined by $y^3 - x^4 = 0$. The associated value semigroup is
\[
\Gamma = \langle 3,4 \rangle = \{0, 3, 4, 6, 7, 8, 9, \dots\},
\]
which contains all integers $n \geq 6$.

Consider the subsemimodule $\Delta = (3,4) = \Gamma \setminus \{0\} \in \mathscr{D}_1$, which is the root of the tree $G_\Gamma$. Applying the deletion maps $d_{1,i}$ corresponding to each minimal generator, we obtain:
\begin{itemize}
    \item $d_{1,1}((3,4)_{\Gamma} ) = (3,4)_{\Gamma}  \setminus \{3\} = \{4,6,7,8,9,\dots\} = (4,6)_{\Gamma}  \in \mathscr{D}_2$,
    \item $d_{1,2}((3,4)_{\Gamma} ) = (3,4)_{\Gamma}  \setminus \{4\} = \{3,6,7,8,9,\dots\} = (3,8)_{\Gamma}  \in \mathscr{D}_2$.
    \item $d_{2,1}((4,6)_{\Gamma} ) = (4,6)_{\Gamma}  \setminus \{4\} = \{6,7,8,9,\dots\} = (6,7,8)_{\Gamma}  \in \mathscr{D}_3. $
\end{itemize}
Note that both $(4,6)$ and $(3,8)$ are distinct elements in $\mathscr{D}_2$, corresponding to two different branches from the root.
\begin{figure}[h]
\centering
\begin{tikzpicture}[scale=0.8, level distance=0.1cm,
   level 1/.style={sibling distance=5cm, level distance=1.5cm},
   level 2/.style={sibling distance=6cm, level distance=2cm},
    level 3/.style={sibling distance=3cm, level distance=2cm},
    level 4/.style={sibling distance=2.3cm, level distance=2cm},
    level 5/.style={sibling distance=2cm, level distance=2cm},
   grow'=up]
 
\node {$(3,4)_{\Gamma} $} 
   child {node {$(4,6)_{\Gamma} $} 
    child {node {$(6,7,8)_{\Gamma} $}  
        child {node {$(7,8,9)_{\Gamma} $}
            child {node {$(8,9,10)_{\Gamma} $}
                child{node {$(9,10,11)_{\Gamma} $}edge from parent node[left]{$d_{1}$}}
                    child {node {$(8,10)_{\Gamma} $}edge from parent node[right]{$d_{2}$}}
                        child {node {$(8,9)_{\Gamma} $}edge from parent node[right]{$d_{3}$}}
                        edge from parent node[left]{$d_{1}$}}
                            child {node {$(7,9)_{\Gamma} $}
                            child [missing]{}
                            child [missing]{}
                                child {node {$(7,12)_{\Gamma} $}edge from parent node[right]{$d_{2}$}}
                                edge from parent node[right]{$d_{2}$}}
                                    child {node {$(7,8)_{\Gamma} $}edge from parent node[right]{$d_{3}$}}
                                    edge from parent node[left]{$d_{1}$}}  
                                        child {node {$(6,8)_{\Gamma} $}
                                            child[missing]{}
                                                child[missing]{}
                                                    child {node {$(6,11)_{\Gamma} $}
                                                    child {node {$(6)_{\Gamma} $}edge from parent node[right]{$d_{2}$}}edge from parent node[right]{$d_{2}$}}
                                                    edge from parent node[right]{$d_{2}$}}   
                                                        child {node {$(6,7)_{\Gamma} $}edge from parent node[right]{$d_{3}$}}edge from parent node[left]{$d_{1}$}}
                                                            child {node {$(4,9)_{\Gamma} $ }
                                                                child {node {$(4)_{\Gamma} $}edge from parent node[right]{$d_{2}$}}edge from parent node[right]{$d_{2}$}}
                                                                    edge from parent node[left]{$d_{1}$}}
                                                                        child {node {$(3,8)_{\Gamma} $} 
                                                                            child {node {$(3)_{\Gamma} $}
                                                                            edge from parent node[right]{$d_{2}$}}edge from parent node[right]{$d_{2}$}
};

\end{tikzpicture}
    \caption{Tree for the case of  $E_{6}$ type singularity}
    \label{fig:enter-label}
\end{figure}
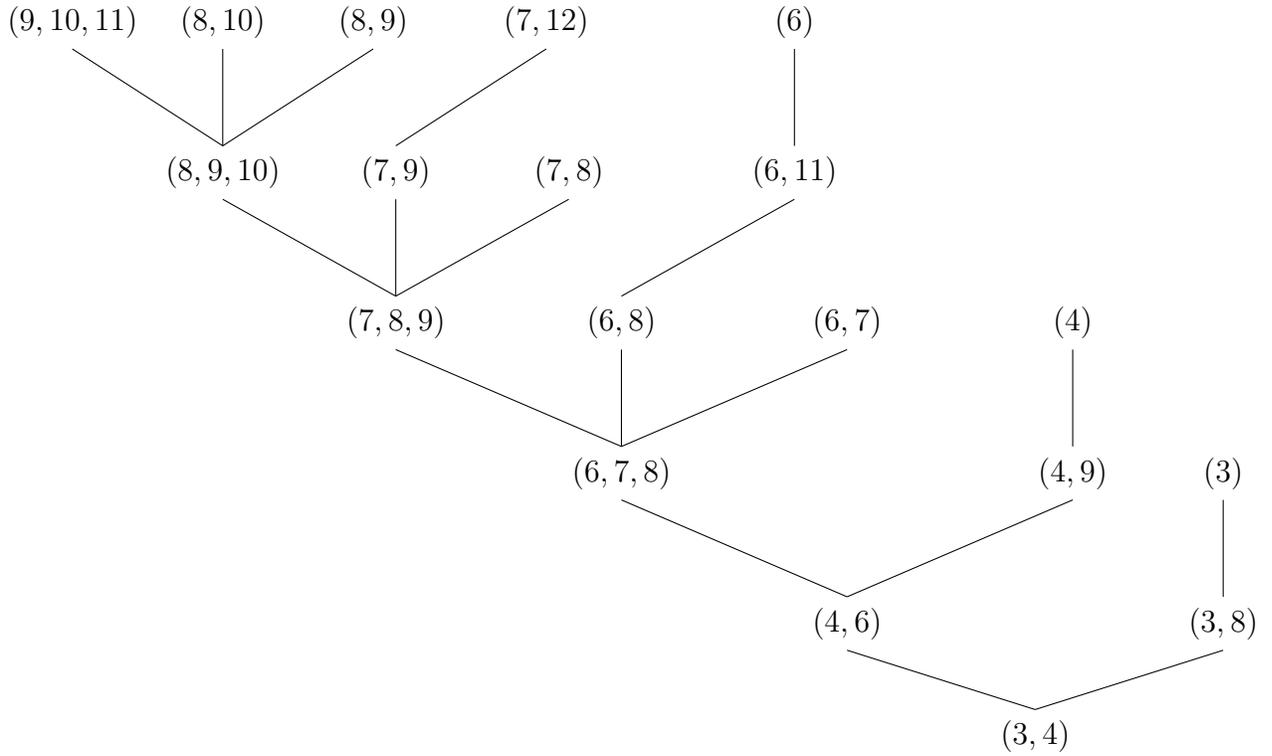\\
\end{example}

\section{Piecewise fibrations induced by the edges of the tree}\label{section: Piecewise fibrations}
In this section, we revisit key results from \cite{hajli2025geometry}. Let $(C,O)$ be the germ of a unibranch plane curve singularity with semigroup $\Gamma$, where $\Gamma$ is either monomial or of the form $\langle \alpha_1, \alpha_2 \rangle$. Write $\Gamma = \langle \alpha_1, \dots, \alpha_e \rangle$. As explained in Section~\ref{section: Tree structure associated with valuation semigroup}, one can then construct a tree from $\Gamma$. 

In \cite{hajli2025geometry}, the authors proved that an edge of the $\Gamma$-subsemimodule tree $G_\Gamma$, connecting a semimodule $\Delta$ to $m(\Delta)$, induces a piecewise trivial fibration between the corresponding strata $C^{[\Delta]}$ and $C^{[m(\Delta)]}$ inside the punctual Hilbert scheme. This fibration reveals a recursive geometric organization of the Hilbert scheme.

Let $\Delta = (\gamma_1,\dots,\gamma_m)_{\Gamma}$. The \emph{syzygy module} of $\Delta$ is the $\Gamma$-subsemimodule defined by
\[
\operatorname{Syz}(\Delta) := \{\sigma \in \Delta \mid \sigma = \gamma_{i_1} + b_1 = \gamma_{i_2} + b_2,\ \exists\, \gamma_{i_1} \neq \gamma_{i_2},\ b_1,b_2 \in \Gamma \},
\]
which is generated as a $\Gamma$-semimodule by a finite set $\{\sigma_1,\dots,\sigma_n\}$. Let $\gamma_\Delta = \max(\Gamma \setminus \Delta)$ and $m(\Delta) = \Delta \cup \{\gamma_\Delta\}$. For any ideal $I \in C^{[\Delta]}$ there exists a generating set of the form
\[
\{f_{\gamma_1}(t),\dots,f_{\gamma_m}(t)\},
\]
where $f_{\gamma_j}(t) = t^{\gamma_j} + \sum_{k \in \Gamma_{>\gamma_j} \setminus \Delta} \lambda_j^{k-\gamma_j} t^{k}$ for some $\lambda_j^{k-\gamma_j} \in \mathbb{C}$.

This yields a canonical morphism
\begin{equation}\label{Piecewise fibration}
C^{[\Delta]} \longrightarrow C^{[m(\Delta)]},\qquad 
\langle f_{\gamma_1}(t),\dots,f_{\gamma_m}(t) \rangle \longmapsto 
\langle f_{\gamma_1}(t),\dots,f_{\gamma_m}(t), t^{\gamma_\Delta} \rangle .
\end{equation}

\begin{theorem}\cite{hajli2025geometry}
Let $C$ be either a plane curve singularity defined by $x^p =y^q $ with $\gcd(p,q) = 1$, 
 or a curve singularity with a monomial semigroup $\Gamma$. Let $\Delta$ be a $\Gamma$-subsemimodule of $\Gamma$, and let $\gamma_\Delta = \max(\Gamma \setminus \Delta)$. Then there exists a canonical morphism
\[
\pi_{\Delta} \colon C^{[\Delta]} \longrightarrow C^{[m(\Delta)]},
\]
which is isomorphic to a trivial fibration over its image. The fiber of $\pi$ is isomorphic to an affine space $\mathbb{A}^{B(\Delta)}$, where
\[
B(\Delta) = \# \left\{ \gamma_i \mid \gamma_i < \gamma_\Delta \right\}
\]
i.e. is number of minimal generators below $\gamma_\Delta$,
and the $\gamma_i$ range over the minimal generators of $\Delta$ as a $\Gamma$-subsemimodule.

\label{piecewise fibration}
\end{theorem}

\begin{corollary}\cite{hajli2025geometry}
Let $C$ be either a plane curve singularity defined by $x^p = y^q $ with $\gcd(p,q) = 1$, 
 or a curve singularity with a monomial semigroup $\Gamma$.   Then given a subsemimodule $\Delta=(\gamma_{1},\dots,\gamma_{m})_{\Gamma}\subset \Gamma$, $Syz(\Delta)=(\sigma_{1},\dots,\sigma_{n})_{\Gamma}$,  we have 
$$[C^{[\Delta]}]=[C^{[m(\Delta)]}]\cdot \mathbb{L}^{|\gamma_{i}|_{\gamma_{i}<\gamma_{\Delta}}-|\sigma_{i}|_{\sigma_{i}<\gamma_{\Delta}}}$$
\label{mononial case,relation in grothendieck}
\end{corollary}

\begin{example}
$C=\{y^{4}=x^{7}\}\subset\mathbb{C}^{2}$,
$\Gamma=\langle 4,7 \rangle= \{4,7,8,11,12,14,15,16,18,\rightarrow\}$, $\Delta=(8,11)_{\Gamma}=\{8,11,12,15,16,17,18,19,20,22,\rightarrow\}$, $\gamma_{\Delta}=21$,  $m(\Delta)=(8,11,21)_{\Gamma}=\{8,11,12,15, 16 ,18\rightarrow\}$. $Syz(\Delta)=(15,32)$. \\

\noindent
Since $\{\gamma_{i}\mid \gamma_{i}<\gamma_{\Delta}\}=\{8,11\}$, $\{\sigma_{i}\mid \sigma_{i}<\gamma_{\Delta}\}=\{15\}$, then 
$[C^{[(8,11)_{\Gamma}]}]=[C^{[(8,11,21)_{\Gamma}]}]\cdot \mathbb{L}$. 
\end{example}

\section{Algorithm for $Z_{(C, O)}^{\text {Hilb}}$} \label{Algorithm 1}
In this section, we introduce Algorithm~\ref{algorithm motivic} for computing the motivic Hilbert zeta function for curve singularities with a monomial valuation semigroup or for singular curves defined by $y^{k}=x^{n}$, where $\gcd(k,n)=1$. The algorithm relies on several subroutines (Algorithms~\ref{algorithm min generators}--\ref{algorithm tree}), each performing a critical step in the overall computation.

We begin by recalling the definition of the Grothendieck ring of varieties and of the motivic Hilbert zeta function.

\begin{definition}
The \emph{Grothendieck ring} of complex varieties, denoted $K_{0}(\mathrm{Var}_{\mathbb{C}})$, is the ring generated by isomorphism classes $[X]$ of complex algebraic varieties $X$, subject to the \emph{scissor relation} (also called the \emph{cut-and-paste relation}):
\[
[X] = [X \setminus Y] + [Y]
\]
for every closed subvariety $Y \subseteq X$. Multiplication is defined by $[X] \cdot [Y] = [X \times Y]$.
\end{definition}

\begin{definition}
Let $C^{[\ell]}$ denote the $\ell$-th punctual Hilbert scheme of the plane curve singularity $(C,O)$, parametrizing ideals of colength $\ell$ supported at the origin. The \emph{motivic Hilbert zeta function} is defined as the generating series
\[
Z_{(C,O)}^{\mathrm{Hilb}}(q) := 1 + \sum_{\ell=1}^{\infty} \bigl[C^{[\ell]}\bigr] \, q^{\ell} \in K_{0}(\mathrm{Var}_{\mathbb{C}})[[q]].
\]
\end{definition}

The stratification of the punctual Hilbert scheme induced by the semimodule filtration yields
\[
\bigl[C^{[\ell]}\bigr] = \sum_{\Delta \in \mathscr{D}_{\ell}} \bigl[C^{[\Delta]}\bigr] \qquad \text{in } K_{0}(\mathrm{Var}_{\mathbb{C}}).
\]
Consequently, the motivic Hilbert zeta function can be rewritten as
\[
Z_{(C,O)}^{\mathrm{Hilb}}(q) = \sum_{\ell \geq 0} \sum_{\Delta \in \mathscr{D}_{\ell}} \bigl[C^{[\Delta]}\bigr] \, q^{\ell}.
\]

Algorithm~\ref{algorithm motivic} takes as input the germ of a unibranch curve singularity $(C,O)$ whose complete local ring is monomial, i.e.\ $\widehat{\mathcal{O}_{C,O}} = \mathbb{C}[[t^{\alpha_{1}},\dots,t^{\alpha_{e}}]]$, and whose valuation semigroup $\Gamma$ is either monomial or of the form $\Gamma = \langle \alpha_{1},\alpha_{2}\rangle$, together with its conductor $c$.
Since $\Gamma$ is infinite, direct computational handling is impractical. We therefore work with the truncated set
\[
\Gamma_{\leq n^{*}} := \Gamma \cap [0, n^{*}],
\]
where $[a,b] := \{ n \in \mathbb{Z} \mid a \leq n \leq b \}$, for a suitable bound $n^{*}$. This truncation makes Algorithms~\ref{algorithm motivic} and~\ref{algorithm tree} executable within practical limits.

Define
\[
N := \min\bigl\{ n^{*} \mid \text{the truncated set } \Gamma_{\leq n^{*}} \text{ suffices for Algorithm~\ref{algorithm motivic}} \bigr\}.
\]
By Theorem~\ref{limite of effective bound}, the choice $n^{*}\geq  (c-1)(\alpha_{1}+2)$ guarantees the correctness of the algorithm.

The output of Algorithm~\ref{algorithm motivic} is the truncation of the motivic Hilbert zeta function up to conductor $c$:
\[
Z_{(C,0)}^{\mathrm{Hilb},\leq c}(q) = \sum_{\Delta \in \mathscr{D}_{1} \cup \dots \cup \mathscr{D}_{c}} \bigl[C^{[\Delta]}\bigr] \, q^{\ell(\Delta)},
\]
where $\ell(\Delta) = \#(\Gamma \setminus \Delta)$.

We initialise $\mathscr{D}_{1} = \{\Delta_{1}\}$ with $\Delta_{1} := \Gamma \setminus \{0\} = (\alpha_{1},\dots,\alpha_{e})$. As explained in Section~\ref{section: Tree structure associated with valuation semigroup}, there is a tree $G_{\Gamma}$ associated with $\Gamma$; using Algorithm~\ref{algorithm tree} (line~1 of Algorithm~\ref{algorithm motivic}), we obtain the sets $\mathscr{D}_{2},\dots,\mathscr{D}_{c}$.

The algorithm then loops over $\ell = 1,\dots,c$, computing each coefficient
\[
\bigl[C^{[\ell]}\bigr] = \sum_{\Delta \in \mathscr{D}_{\ell}} \bigl[C^{[\Delta]}\bigr]
\]
of the motivic Hilbert zeta function.

Inside the main loop, for each $\Delta \in \mathscr{D}_{\ell}$ we perform the following steps:
\begin{enumerate}
\item Compute $\gamma_{\Delta} := \max(\Gamma \setminus \Delta)$ (the Frobenius element of $\Delta$).
\item Use Algorithm~\ref{algorithm syz} to find $\operatorname{Syz}(\Delta)$, the syzygy module of $\Delta$.
\item Using Algorithm~\ref{algorithm min generators}, determine
      \begin{itemize}
        \item $T_{\Delta} = \{\gamma_{1},\dots,\gamma_{m}\}$, the minimal generators of $\Delta$;
        \item $T_{\operatorname{Syz}(\Delta)} = \{\sigma_{1},\dots,\sigma_{n}\}$, the minimal generators of $\operatorname{Syz}(\Delta)$.
      \end{itemize}
\item For the contribution counts, compute
      \[
      C_{1} = \#\{\gamma_i \mid \gamma_i < \gamma_{\Delta}\}, \qquad
      C_{2} = \#\{\sigma_i \mid \sigma_i < \gamma_{\Delta}\}.
      \]
\end{enumerate}

\begin{algorithm}
\caption{Motivic Hilbert zeta function of curve singularity}
\begin{algorithmic}[1]
\REQUIRE Let $C$ be a curve singularity with valuation semigroup $\Gamma=\langle\alpha_{1},\dots,\alpha_{e}\rangle$. Let $c$ be the conductor of $\Gamma$.                  
\ENSURE Motivic Hilbert zeta function \[Z_{(C,0)}^{Hilb, \leq c}(q)=\sum_{\Delta\in \mathscr{D}_{1},\dots, \mathscr{D}_{c}}[C^{[\Delta]}]q^{\ell(\Delta)}.\]     
 \STATE   $\mathscr{D}_{1}=\{\Delta_{1}=(\alpha_{1},\dots,\alpha_{e})\}$ , $\mathscr{D}_{2}=\{\Delta_{2,i}=(\Delta_{1}\setminus\{\alpha_{i})\}\},\dots,\mathscr{D}_{c};$\\ 

\FOR {$ \ell=1,\dots c$}
\FOR{$\Delta\in \mathscr{D}_{\ell}$} 
\STATE
$\gamma_{\Delta}=max(\Gamma\setminus \Delta)\leftarrow{\mbox{Frobenius element of } \Delta};$ 
\STATE
$T_{\Delta}=\{\gamma_{1},\dots,\gamma_{m}\}\leftarrow{\mbox{minimal generators of  } \Delta};$
\STATE
compute $Syz(\Delta);$
\STATE
$T_{Syz(\Delta)}=\{\sigma_{1},\dots,\sigma_{n}\}\leftarrow{\mbox{minimal generators of  } Syz(\Delta)};$
\STATE
$m(\Delta)=\Delta\cup \{\gamma_{\Delta}\};$ 
\STATE
$C_1 = 0;$
\FOR{$\gamma_{i}\in T_{\Delta}$}
\IF{$\gamma_{i}<\gamma_{\Delta}$}
\STATE
$C_1 \leftarrow C_1 + 1;$
\ENDIF
\ENDFOR
\STATE
$|\gamma_{i}|_{\gamma_{i}<\gamma_{\Delta}}\leftarrow C_{1};$ 
\STATE
$C_2 = 0;$
\FOR{$\sigma_{i}\in T_{Syz(\Delta)}$}
\IF{$\sigma_{i}<\gamma_{\Delta}$}
\STATE
$C_{2}\leftarrow C_{2}+1;$
\ENDIF
\ENDFOR
\STATE
$|\sigma_{i}|_{\sigma_{i}<\gamma_{\Delta},\sigma_{i}\prec\gamma_{\Delta}}\leftarrow C_{2};$
\STATE
$[C^{[\Delta]}]\leftarrow [C^{[m(\Delta)]}]\cdot \mathbb{L}^{|\gamma_{i}|_{\gamma_{i}<\gamma_{\Delta}}-|\sigma_{i}|_{\sigma_{i}<\gamma_{\Delta},\sigma_{i}\prec\gamma_{\Delta}}};$
\ENDFOR
\STATE
$[C^{\ell}]=\sum_{\Delta\in \mathscr{D}_{\ell}}[C^{\Delta}];$
\ENDFOR
\STATE
${Z_{(C,0)}^{Hilb}(q)}\leftarrow 1+\sum_{\ell\geq 1}[C^{[\ell]}]q^{\ell};$
\RETURN $Z_{(C,0)}^{Hilb,\leq c}(q).$
\end{algorithmic}
\label{algorithm motivic}
\end{algorithm}

\begin{algorithm}[h]
\caption{ Minimal generators of $\Delta$ as $\Gamma$-module, $T_{\Delta}$}
\begin{algorithmic}[1]
\REQUIRE Let $C$ be a curve singularity with valuation semigroup $\Gamma=\langle\alpha_{1},\dots,\alpha_{e}\rangle$. Let $c$ be the conductor of $\Gamma$, $\Delta$ be a  subsemimodule of $\Gamma$. We assume that $\Delta=\{a_{1},\dots,a_{s}\}$ such that $a_{i}<a_{i+1}$,  where $s=\#(\Delta)$ and $\#(\Gamma\setminus \Delta)=\ell$. 
\ENSURE Minimal generators of $\Delta$ as $\Gamma$-module $T_{\Delta}$.\\      
 \STATE let $T_{\Delta}=\{a_{1}\}$;\\ 
\FOR {$  i=2,...,s$}
\FOR{$x \in T_{\Delta}$} 

\STATE
$r_{i} \leftarrow a_{i} - x$;
\IF {$r_{i}\in \Gamma$:}
\STATE
    break; 
\ENDIF
\IF {$r_{i}\notin \Gamma$}
\STATE
$T_{\Delta}\leftarrow  T_{\Delta}\cup \{a_{i}\}$;
\ENDIF
\ENDFOR

\ENDFOR

\RETURN $T_{\Delta}.$
\end{algorithmic}
\label{algorithm min generators}
\end{algorithm}

\begin{algorithm}
\caption{Syzygy of $\Delta$, $Syz(\Delta)$}
\begin{algorithmic}[1]
\REQUIRE Let $C$ be a curve singularity with valuation  semigroup $\Gamma=\langle\alpha_{1},\dots,\alpha_{e}\rangle$. Let $c$ be the conductor of $\Gamma$, $\Delta$ be a subsemimodule of $\Gamma$. We assume that $\Delta=\{a_{1},\dots,a_{s-1}\}$ such that $a_{i}<a_{i+1}$,  where $s=\#(\Delta)$.  Minimal generators of $\Delta$ as $\Gamma$-module, $T_{\Delta}=\{\gamma_{1},\dots,\gamma_{m}\}$. 
\ENSURE  Syzygy of $\Delta$, $Syz(\Delta)$. \\      
 \STATE let $Syz(\Delta)=\{\}$;\\ 
\FOR {$i=1,...,s-1$}
\STATE
$count = 0$;
\FOR{$j = 0,\dots, m-1$} 

\STATE
$r_{i,j} \leftarrow a_{i} - \gamma_{j}$;
\IF {$r_{i,j}\in \Gamma$}
\STATE
    $count \leftarrow count + 1$;
\ENDIF

\IF {$count >=2$}
\STATE
break;
\ENDIF

\ENDFOR
\IF {$count = 2$}
\STATE
$Syz(\Delta)\leftarrow Syz(\Delta)\cup \{a_{i}\};$
\ENDIF
\ENDFOR

\RETURN $Syz(\Delta).$
\end{algorithmic}
\label{algorithm syz}
\end{algorithm}
\clearpage
\begin{algorithm}
\caption{ The tree obtained from $\Gamma$, $G_{\Gamma}$}
\begin{algorithmic}[1]
\REQUIRE Let $C$ be a curve singularity with valuation semigroup $\Gamma=\langle\alpha_{1},\dots,\alpha_{e}\rangle$. Let $c$ be the conductor of $\Gamma$. 
\ENSURE $\ell = 1,\dots,c$, $\mathscr{D}_{\ell}$, the set of $\Gamma$-subsemimodule $\Delta$ such that $\#(\Gamma\setminus \Delta)=\ell$.  \\
 \STATE let $\mathscr{D}_{1}=\{\Delta_{1}=(\alpha_{1},\dots,\alpha_{e})\}$, $(\alpha_{1},\dots,\alpha_{e})=\Gamma\setminus \{0\}$;\\ 
\FOR {$  \ell=1,\dots,c-1$}
\STATE
$\mathscr{D}_{\ell+1}=\{\};$
\FOR{$\Delta\in \mathscr{D}_{\ell}$} 
\STATE
compute the set of minimal generators of $\Delta$;  $T_{\Delta}=\{\gamma_{1},\dots,\gamma_{m}\}$, $m\leq \alpha_{1}$; 
\FOR{$\gamma_{i}\in T_{\Delta}$}
\STATE
$\Delta_{i}=\Delta \setminus \{\gamma_{i}\};$
\STATE
$\mathscr{D}_{\ell+1}\leftarrow \mathscr{D}_{\ell+1}\cup \{\Delta_{i}\};$
\ENDFOR
\ENDFOR
\ENDFOR

\RETURN $\mathscr{D}_{\ell},\ell = 1,\dots,c$.
\end{algorithmic}
\label{algorithm tree}
\end{algorithm}

In the final part of this section, we analyze the time complexity of the algorithms.

\begin{lemma}
\label{complexity min generator}
Let $\Delta$ be a $\Gamma$-subsemimodule with $\#(\Gamma\setminus\Delta)=\ell$.  
The time complexity of Algorithm~\ref{algorithm min generators} is $O\bigl(n^{*}\alpha_{1}(n^{*}-\ell)\bigr)$. 
\begin{proof}
The outer loop (for $i = 2,\dots, s$ where $s = n^{*}-\ell$) runs at most $O(n^{*} - \ell)$ times. 
For the inner loop (over $x\in T_{\Delta}$), the size of the set of minimal generators of the $\Gamma$-module is at most $\alpha_{1}$. 
Each iteration of the inner loop requires a membership check in $\Gamma$, which takes $O(n^{*})$ time. 
Thus, the overall time complexity is $O\bigl(n^{*}\alpha_{1}(n^{*}-\ell)\bigr)$. 
\end{proof}
\end{lemma}    

For the same reason, the time complexity of Algorithm~\ref{algorithm syz} is $O\bigl(n^{*}\alpha_{1}(n^{*}-\ell)\bigr)$. Let $\#\mathscr{D}_{\ell} = \chi_{\ell}$. 

\begin{lemma}\label{complexity tree}
The time complexity of Algorithm~\ref{algorithm tree} is 
\[
\sum_{\ell=1}^{c-1} O\bigl(\chi_{\ell} \cdot n^{*}\alpha_{1}(n^{*}-\ell)\bigr).
\]
\begin{proof}
One only needs to notice that, on line~7, deleting an element from $\Delta$ costs $O(1)$ time. 
The reason is that if $\Delta$ is treated as a sorted list, the element to be deleted (being one of the generators of $\Delta$ as a $\Gamma$-module) must appear near the beginning of the list, so removal takes $O(1)$ time. 
If $\Delta$ is treated as a set, deleting a specified element also takes $O(1)$ time. 
Hence the total complexity is
\[
\sum_{\ell=1}^{c-1}\bigl[O(\chi_{\ell} \cdot n^{*}\alpha_{1}(n^{*}-\ell)) + O(\chi_{\ell} \cdot \alpha_{1})\bigr]
= \sum_{\ell=1}^{c-1} O\bigl(\chi_{\ell} \cdot n^{*}\alpha_{1}(n^{*}-\ell)\bigr).
\]
\end{proof}
\end{lemma}

\begin{theorem}
The time complexity of Algorithm~\ref{algorithm motivic} is 
\[
\sum_{\ell=1}^{c} O\bigl(c \cdot \chi_{\ell} \cdot n^{*} \cdot \alpha_{1} \cdot (n^{*}-\ell)\bigr).
\]
\begin{proof}
By Lemma~\ref{complexity min generator} and Lemma~\ref{complexity tree}, the overall time complexity is
\[
\sum_{\ell=1}^{c} O\bigl(c \cdot \chi_{\ell} \cdot n^{*} \cdot \alpha_{1} \cdot (n^{*}-\ell)\bigr).
\]
\end{proof}
\end{theorem}

\section{Algorithm for $Zm_{(C, O)}^{\text {Hilb}}$}\label{section Algorithm for Zm}

In this section, we  introduce Algorithm~\ref{new algorithm motivic}, Algorithm~\ref{algm: C^Delta,<m}, and Algorithm~\ref{algm: C^Delta,i} to compute the \emph{new} motivic Hilbert zeta function $Zm_{(C, O)}^{\text {Hilb}}$ for curve singularities $C$ defined by $y^{p}=x^{q}$, where $\gcd(p,q)=1$.

We begin by recalling the definition of the new motivic Hilbert zeta function introduced in \cite{hajli2025geometry}, along with some related properties of subvarieties of the punctual Hilbert scheme with a fixed number of generators.

Let $(C,O)$ be the germ of a unibranch curve singularity defined over $\mathbb{C}$, and let $A = \mathcal{O}_{C,O}$ denote its local ring at $O$.  
Let $m(I)$ be the minimal number of generators of ideal $I$. For $m\in \mathbb{Z}_{\geq 1}$, denote the subvariety of punctual Hilbert scheme parametrizing ideals of $A$ with $m$ minimal number of generators:  
\[C^{[\ell],m}:=\{I\in C^{[\ell]}\mid m(I)=m\}.\]

\begin{definition}\cite{hajli2025geometry}
The \emph{motivic Hilbert zeta function with fixed number of generators} is defined as
\begin{align}
Zm_{(C,O)}^{\mathrm{Hilb}}(a^2, q^2) 
&:= \sum_{\ell \geq 0} \sum_{m \geq 1} q^{2\ell} (1 - a^2)^{m-1} [C^{[\ell],m}] \nonumber
 \\
&= \sum_{\ell \geq 0} \sum_{\Delta \in \mathscr{D}_\ell} \sum_{m \geq 1} q^{2\ell} (1 - a^2)^{m-1} [C^{[\Delta],m}],\nonumber
\end{align}
where 
\[
\mathscr{D}_\ell = \big\{ \Delta \subset \Gamma \mid \Delta \text{ is a } \Gamma\text{-semimodule},\ \#(\Gamma \setminus \Delta) = \ell \big\},
\]
and $C^{[\Delta],m}$ denotes the locally closed subvariety of $C^{[\Delta]}$ parametrizing ideals $I \subset \mathcal{O}_C$ with exactly $m$ minimal generators and value semimodule $\Delta$, i.e.,
\[
C^{[\Delta],m} = \{ I \in C^{[\Delta]} \mid m(I) = m \}.
\]
\end{definition}

We now introduce the following results, which will be used to construct the algorithm for the new motivic Hilbert zeta function:

\begin{definition}\cite{hajli2025geometry}\label{defCDeltai}
 Let $\Delta$ be a $\Gamma$-subsemimodule with a minimal generating set $\{\gamma_1,\dots,\gamma_n\}$. 
For $\underline{i} = \{i_1, \dots, i_m\} \subseteq \{1, \dots, n\}$,  define
\[
C^{[\Delta],\underline{i}} := \left\{ I_\lambda \in C^{[\Delta]} \;\middle|\; 
I_\lambda \text{ admits a generating set of size } m \text{ consisting of } f_{\gamma_{i_1}}(\lambda), \dots, f_{\gamma_{i_m}}(\lambda) \right\}.
\]
\end{definition}

\begin{proposition}\label{decomposition of C Delta m}\cite{hajli2025geometry}

\item 
\begin{enumerate}
    \item  For $m \in \mathbb{Z}_{\geq 1}$, we have the stratification
    \[
    C^{[\Delta],m} = C^{[\Delta],\leq m} \setminus C^{[\Delta],\leq m-1},
    \]
    where $C^{[\Delta],\leq m} = \{ I \in C^{[\Delta]} \mid m(I) \leq m \}$.

    \item For each subset $\underline{i} \subseteq \{1, \dots, n\}$, the strata $C^{[\Delta],\underline{i}}$, $C^{[\Delta],\leq m}$, and $C^{[\Delta],m}$ are locally closed subvarieties of $C^{[\Delta]}$.

    \item For subsets $\underline{i}, \underline{j} \subseteq \{1, \dots, n\}$ with $\# \underline{i} = \# \underline{j} = m$, we have
    \[
    C^{[\Delta],\underline{i}} \cap C^{[\Delta],\underline{j}} = C^{[\Delta],\underline{i} \cap \underline{j}}.
    \]
       Furthermore,   denote     $\{ \underline{i} \subseteq \{1,\dots,n\} \mid  \# \underline{i} = m \} = \{ \underline{i}_{1},\dots, \underline{i}_{t}\}$  the set of subsets $\{1,\dots,n\}$ of size $m$ with $t = \binom{n}{m}$.  By the inclusion-exclusion principle, we have

    \begin{align*}
    [C^{[\Delta],\leq m}] 
    &= \left[\bigcup_{\substack{\underline{i} \subseteq \{1,\dots,n\} \\ \#\underline{i} = m}} C^{[\Delta],\underline{i}}\right] \\
    &= \sum_{\substack{\underline{i}}} [C^{[\Delta],\underline{i}}] 
       - \sum_{\substack{\underline{i},\underline{j} \\ \underline{i} < \underline{j}}} [C^{[\Delta],\underline{i} \cap \underline{j}}] 
       + \cdots 
       + (-1)^{t-1} [C^{[\Delta],\underline{i}_1 \cap \cdots \cap \underline{i}_t}] \\
    &= \sum_{\substack{\mathcal{I} \subseteq \{1,\dots,t\} \\ \mathcal{I} \neq \varnothing}} (-1)^{|\mathcal{I}|-1} [C^{[\Delta],\bigcap_{e \in \mathcal{I}} \underline{i}_e}] \quad \text{in } K_0(\mathrm{Var}_{\mathbb{C}}),
    \end{align*}
    where the ordering $\underline{i} < \underline{j}$ is lexicographic.
\end{enumerate}\end{proposition}

\begin{theorem}\label{decomposition of C Delta i}\cite{hajli2025geometry}
Let $\Gamma $  be the value semigroup associated with the germ of a curve singularity defined by $y^p=x^q$ with $\gcd(p,q) = 1$.  Let $\Delta = (\gamma_1, \dots, \gamma_n)_\Gamma$, with $n > 1$, be a $\Gamma$-subsemimodule minimally generated by $n$ elements as a $\Gamma$-module.  Fix an integer $m$ such that $1 < m \leq n$, and let
\[
\underline{i} = \{i_1, \dots, i_m\} \subset \{1, \dots, n\}, \quad \{i'_1, \dots, i'_{n-m}\} = \{1, \dots, n\} \setminus \underline{i}.
\]
For each $e = 1, \dots, n-m$, choose a minimal generator in the $Syzgy$ set of $\Gamma-$ semimodule  $(\gamma_{i_1}, \dots, \gamma_{i_m}, \gamma_{i'_1}, \dots, \gamma_{i'_{e-1}})_{\Gamma}$:
\[
 \sigma_{i_{j_e}} \in T_{\operatorname{Syz}((\gamma_{i_1}, \dots, \gamma_{i_m}, \gamma_{i'_1}, \dots, \gamma_{i'_{e-1}})_{\Gamma})}
\]
such that $\sigma_{i_{j_e}} < \gamma_{i'_e}$. 

Let $Y_{\underline{i}_{j_e}}$ be the subvariety defined by the conditions:
\[
\sum_{k \in \{i_1, \dots, i_m, i_{1}^{\prime},\dots, i_{e-1}^{\prime}\}} (\mathcal{G}^{(e)}_{\lambda})_k \circ (\mathcal{S}^{(e)}_{\nu})_{i_{j_e}}^k \equiv 0 \mod t^{\gamma_{i'_e}}, 
\quad \text{and} \quad 
(Eq^{(e)})^{\gamma_{i'_e} - \sigma_{i_{j_e}}}_{i_{j_e}} \neq 0,
\]
where $(Eq^{(e)})^{\gamma_{i'_e} - \sigma_{i_{j_e}}}_{i_{j_e}}$ is the coefficient of $\phi_{\gamma_{i'_e}}$ in the expansion of the left-hand side, as defined in the proof.

Now define
\[
Y_{\underline{i}_{\underline{j}}} = \bigcap_{e=1}^{n-m} Y_{\underline{i}_{j_e}}.
\]
Then we have the following decomposition:
\[
C^{[\Delta],\leq m} = \bigcup_{\substack{\underline{i} \subseteq \{1,\dots,n\} \\ |\underline{i}| = m}} C^{[\Delta],\underline{i}} = \bigcup_{\substack{\underline{i} \subseteq \{1,\dots,n\} \\ |\underline{i}| = m}} \bigcup_{\underline{j}} Y_{\underline{i}_{\underline{j}}},
\]
where the second union runs over all tuples $\underline{j} = (j_1, \dots, j_{n-m})$ corresponding to valid tuples chosen in the Syzygy  as above.

Furthermore, we have:
\[
Y_{\underline{i}_{\underline{j}}} \cong ({\mathbb{C}}^*)^{n - m} \times_{\operatorname{Spec} {\mathbb{C}}} \mathbb{A}^{N(\Delta) - n + m},
\]
where $N(\Delta) = \dim C^{[\Delta]}$ is the dimension of the stratum, and $C^{[\Delta]} \cong \mathbb{A}^{N(\Delta)}$ as established in \parencite[Theorem 13]{oblomkov2018hilbert}

\end{theorem}

\begin{remark}\cite{hajli2025geometry}\label{inclusion-exclusion principle for C_delta_i}
Assume $C^{[\Delta],\underline{i}}$ is a union of  $Y_{\underline{i}_{\underline{j_{g}}}}$ as described in  Theorem \ref{decomposition of C Delta i}, $\underline{i}_{\underline{j_{g}}}=\underline{i}_{\underline{j_{1}}}, \dots, \underline{i}_{\underline{j_{\eta}}}$. By inclusion-exclusion principle,  we have:
\[ [C^{[\Delta],\underline{i}}]=[\bigcup_{\underline{i}_{\underline{j_{g}}}=\underline{i}_{\underline{j_{1}}},\dots \underline{i}_{\underline{j_{\eta}}}}Y_{\underline{i}_{\underline{j_{g}}}}]= 
\sum_{\underline{j_{g}}}  [Y_{\underline{i}_{\underline{j_{g}}}}]-\sum_{\underline{i}_{\underline{j_{g}}},\underline{i}_{\underline{j_{f}}}, g<f}  [Y_{\underline{i}_{\underline{j_{g}}}}\cap Y_{\underline{i}_{\underline{j_{f}}}}]+\dots + (-1)^{(\eta-1)} [Y_{\underline{i}_{\underline{j_{1}}}}\cap \dots \cap Y_{\underline{i}_{\underline{j_{\eta}}}}]
\]
\end{remark}

In the following, we regard $\sigma_{\underline{i}_{\underline{j_{g}}}}=(\{\sigma_{i_{j_{g_{1}}}}\},\dots,\{\sigma_{i_{j_{g_{n-m}}}}\})$ as $(n-m)$-tuple of sets whose components are sets. We can define a union operator of two tuples of sets: \[\sigma_{\underline{i}_{\underline{j_{g}}}}\cup \sigma_{\underline{i}_{\underline{j_{f}}}}:=(\{\sigma_{i_{j_{g_{1}}}}\}\cup \{\sigma_{i_{j_{f_{1}}}}\} ,\dots,\{\sigma_{i_{j_{g_{n-m}}}}\}\cup \{\sigma_{i_{j_{f_{n-m}}}}\})\]

We define the cardinality of a tuple of sets as the sum of the cardinalities of its components, i.e., 
\[\#(\sigma_{\underline{i}_{\underline{j_{g}}}}\cup \sigma_{\underline{i}_{\underline{j_{f}}}}):=\sum_{e=1}^{n-m}\#(\{\sigma_{i_{j_{g_{e}}}}\}\cup \{\sigma_{i_{j_{f_{e}}}}\}). \]

\begin{corollary}\label{intersection_for Y}\cite{hajli2025geometry}
For the curve $C$ defined by $y^p = x^q$, we assume that $C^{[\Delta],\underline{i}}$ is a union of the sets $Y_{\underline{i}_{\underline{j_g}}}$ as described in Theorem~\ref{decomposition of C Delta i}, where $\underline{i}_{\underline{j_g}} = \underline{i}_{\underline{j_1}}, \dots, \underline{i}_{\underline{j_\eta}}$. Then, for $s \leq \eta$, we have:
\[
Y_{\underline{i}_{\underline{j}_{1}}} \cap  \dots \cap  Y_{\underline{i}_{\underline{j}_{s}}} \cong (\mathbb{C}^{*})^{\#(\sigma_{\underline{i}_{\underline{j}_{1}}} \cup \dots \cup \sigma_{\underline{i}_{\underline{j}_{s}}})} \times _{\mathbb{C}} \mathbb{A}^{N(\Delta) - \#(\sigma_{\underline{i}_{\underline{j}_{1}}} \cup \dots \cup \sigma_{\underline{i}_{\underline{j}_{s}}})}.
\]
\end{corollary}

\begin{algorithm}
\caption{New motivic Hilbert zeta function of curve defined by $y^p=x^q$}
\begin{algorithmic}[1]
\REQUIRE Let $C$ be a curve singularity defined by $y^p=x^q$  with valuation semigroup $\Gamma=\langle p, q \rangle$ with the conductor $c=(p-1)(q-1)$.
\ENSURE New motivic Hilbert zeta function \[Zm_{(C,0)}^{\text{Hilb},\leq c}(a^2,q^2)= \sum_{\ell \geq 0}^{c} \sum_{\Delta \in \mathscr{D}_\ell} \sum_{m \geq 1} q^{2\ell} (1 - a^2)^{m-1} [C^{[\Delta],m}].\]

\STATE $\mathscr{D}_{1}=\{\Delta_{1}=(p,q)_{\Gamma}\}$, $\mathscr{D}_{2}=\{\Delta_{2,1}=(\Delta_{1}\setminus\{p\}),\Delta_{2,2}=(\Delta_{1}\setminus\{q\}) \}$, $\dots$, $\mathscr{D}_{c}$;
\STATE $Zm_{(C,O)}^{\text{Hilb}}(a^2,q^2) = 0 \leftarrow \mbox{output}$;

\FOR{$\ell=1$ \TO $c$}
    \FOR{$\Delta\in \mathscr{D}_{\ell}$} 
        \STATE $T_{\Delta}=\{\gamma_{1},\dots,\gamma_{n}\} \leftarrow{\mbox{minimal generators of } \Delta}$;
        \STATE compute $Syz(\Delta)$;
        \STATE $T_{Syz(\Delta)}=\{\sigma_{1},\dots,\sigma_{n}\} \leftarrow{\mbox{minimal generators of } Syz(\Delta)}$;
        \STATE compute $N(\Delta)=\sum_{i=1}^n \#(\Gamma_{>\gamma_{i}}\setminus \Delta)-\sum_{i=1}^n \#(\Gamma_{>\sigma_{i}}\setminus \Delta)$;
        \STATE $Z_{\text{temp}} =  0$;
        \FOR{$\frac{n }{2} < m \leq n $} 
        	\STATE 
	compute $[C^{[\Delta],\leq m}]$ and $[C^{[\Delta],\leq m-1}]$  $\leftarrow $ Algorithm 6;
            \STATE $Z_{\text{temp}} \leftarrow  Z_{\text{temp}} +  q ^{2\ell}  (1-a^2)^{m-1}([C^{[\Delta],\leq m}]-[C^{[\Delta],\leq m-1}])$; 
        \ENDFOR
        \STATE 
        $Zm_{(C,0)}^{\text{Hilb}}(a^2,q^2) \leftarrow Zm_{(C,0)}^{\text{Hilb}}(a^2,q^2) + Z_{\text{temp}}$;
    \ENDFOR
\ENDFOR

\RETURN $Zm_{(C,0)}^{\text{Hilb},\leq c}(a^2,q^2).$
\end{algorithmic}
\label{new algorithm motivic}
\end{algorithm}

\begin{algorithm}\label{algm:C_Delta_small than m}
\caption{ $[C^{[\Delta],\leq m}]$ }
\begin{algorithmic}[1]
\REQUIRE Let $C$ be a curve singularity defined by $y^p=x^q$   with valuation semigroup $\Gamma=\langle p, q \rangle$. Let $\Delta$  be a $\Gamma$-subsemimodule which is minimally generated by $\gamma_1, \dots, \gamma_n$. Let $Syz(\Delta)$ be the Syzygy set of $\Delta$.  Denote the dimension of $C^{[\Delta]}$ by $N(\Delta)=\sum_{i=1}^n \#(\Gamma_{>\gamma_{i}}\setminus \Delta)-\sum_{i=1}^n \#(\Gamma_{>\sigma_{i}}\setminus \Delta)$. Let $m\in \mathbb{Z}_{\geq 1}$.  Let  $\{ \underline{i} \subseteq \{1,\dots,n\} \mid  \# \underline{i} = m \} = \{ \underline{i}_{1},\dots \underline{i}_{t}\}$  the set of subsets $\{1,\dots,n\}$ of size $m$ with $t = \binom{n}{m}$.

\ENSURE  $[C^{[\Delta],\leq m}].$ \\      
\STATE
$[C^{[\Delta],\leq m}]= 0$;
\FOR{$ \mathcal{I} \subseteq \{1,\dots, t\}$ }
	\STATE $[C^{[\Delta],\bigcap_{e \in \mathcal{I}} \underline{i}_e}]$ $\leftarrow $ Algorithm 7;
   \STATE  $[C^{[\Delta],\leq m}] \leftarrow [C^{[\Delta],\leq m}] +  (-1)^{|\mathcal{I}|-1}  [C^{[\Delta],\bigcap_{e \in \mathcal{I}} \underline{i}_e}]$;
\ENDFOR

\RETURN $[C^{[\Delta],\leq m}].$
\end{algorithmic}
\label{algm: C^Delta,<m}
\end{algorithm}

\begin{algorithm}
\caption{ $[C^{[\Delta], \underline{i}}]$ }
\begin{algorithmic}[1]
\REQUIRE  
Let $C$ be a curve singularity defined by $y^p=x^q$   with valuation semigroup $\Gamma=\langle p, q \rangle$. Let $\Delta$  be a $\Gamma$-subsemimodule which is minimally generated by $\gamma_1, \dots, \gamma_n$. Let $Syz(\Delta)$ be the syzygy of $\Delta$.  Denote the dimension of $C^{[\Delta]}$ by $N(\Delta)=\sum_{i=1}^n \#(\Gamma_{>\gamma_{i}}\setminus \Delta)-\sum_{i=1}^n \#(\Gamma_{>\sigma_{i}}\setminus \Delta)$. Let $m\in \mathbb{Z}_{\geq 1}$ and $\underline{i}=\{i_{1},\dots,i_{m}\}\subset \{1,\dots n\}$. 

\ENSURE  $[C^{[\Delta], \underline{i}}]$. \\      
\STATE
$[C^{[\Delta], \underline{i}}]= 0$;

\STATE  compute $ \{\gamma_{i_{1}}^{\prime}, \dots, \gamma_{i_{n-m}}^{\prime} \}\leftarrow \{\gamma_1,\dots, \gamma_n\} \setminus \{\gamma_{i_{1}} ,\dots,\gamma_{i_{m}} \}$;

\STATE
Let $\{ \{ \sigma_{i_{j_{1}}}       ,\dots,      \sigma_{i_{j_{n-m}}}  \}\subset T_{Syz(\Delta)} \mid  \sigma_{i_{j_{e}}}  <\gamma_{i_{e}}^{\prime} \}= \{\sigma_{\underline{i}_{\underline{j_1}}}, \dots, \sigma_{\underline{i}_{\underline{j_\eta}}} \}$;  

\FOR{ $s\leq \eta$}
\STATE
\[
[Y_{\underline{i}_{\underline{j}_{1}}} \cap  \dots \cap  Y_{\underline{i}_{\underline{j}_{s}}}] = (\mathbb{L}-1)^{\#(\sigma_{\underline{i}_{\underline{j}_{1}}} \cup \dots \cup \sigma_{\underline{i}_{\underline{j}_{s}}})} \cdot  \mathbb{L}^{N(\Delta) - \#(\sigma_{\underline{i}_{\underline{j}_{1}}} \cup \dots \cup \sigma_{\underline{i}_{\underline{j}_{s}}})};
\]
\ENDFOR
\STATE
\[ [C^{[\Delta],\underline{i}}]\leftarrow [\bigcup_{\underline{i}_{\underline{j_{g}}}=\underline{i}_{\underline{j_{1}}},\dots \underline{i}_{\underline{j_{\eta}}}}Y_{\underline{i}_{\underline{j_{g}}}}]= 
\sum_{\underline{j_{g}}}  [Y_{\underline{i}_{\underline{j_{g}}}}]-\sum_{\underline{i}_{\underline{j_{g}}},\underline{i}_{\underline{j_{f}}}, g<f}  [Y_{\underline{i}_{\underline{j_{g}}}}\cap Y_{\underline{i}_{\underline{j_{f}}}}]+\dots + (-1)^{(\eta-1)} [Y_{\underline{i}_{\underline{j_{1}}}}\cap \dots \cap Y_{\underline{i}_{\underline{j_{\eta}}}}]
;\]

\RETURN $[C^{[\Delta], \underline{i}}].$
\end{algorithmic}
\label{algm: C^Delta,i}
\end{algorithm}
\clearpage

\section{Effective finite length of the value semigroup of a curve singularity}\label{Effective finite length}

Since $\Gamma = \langle \alpha_1, \dots, \alpha_e \rangle$ is a numerical semigroup, we have $\#(\mathbb{N} \setminus \Gamma) < \infty$, but $\Gamma$ itself is infinite. For computational implementation, we face the difficulty that a computer cannot handle an infinite set directly. We therefore approximate $\Gamma$ by truncating it to a finite subset, allowing Algorithm~\ref{algorithm motivic} to operate effectively. Concretely, we choose a positive integer $n^*$ so that the truncated set $\Gamma_{\leq n^*} = \Gamma \cap [0,n^*]$ is sufficient for the algorithm to work correctly.

\begin{definition}
If $N$ is the smallest integer for which Algorithm~\ref{algorithm motivic} works correctly using the truncated set $\Gamma_{\leq N}$, we call $N$ the \emph{minimal effective bound} of $\Gamma$.
\end{definition}

Because Algorithm~\ref{algorithm motivic} depends on Algorithms~\ref{algorithm min generators}, \ref{algorithm syz} and~\ref{algorithm tree}, determining the minimal effective bound $N$ requires a detailed analysis of each of these subroutines.

\begin{lemma}
If the truncated set $\Gamma_{\leq N}$ is sufficient for Algorithm~\ref{algorithm min generators} to operate effectively for every $\Delta \in \mathscr{D}_{\ell}$ with $1 \leq \ell \leq c$, then $N$ can be taken to satisfy
\[
N \le \gamma_{\Delta_c} + c - 1,
\]
where $\gamma_{\Delta_c} = \max\{ \gamma_{\Delta} \mid \Delta \in \mathscr{D}_c \}$ and $\gamma_\Delta = \max(\Gamma \setminus \Delta)$.
\begin{proof}
Choose $\Delta_c \in \mathscr{D}_c$ such that $\gamma_{\Delta_c} = \max\{ \gamma_{\Delta} \mid \Delta \in \mathscr{D}_c \}$. It suffices to show that Algorithm~\ref{algorithm min generators} works for this $\Delta_c$ when we use the truncation $\Gamma_{\leq N'}$ with $N' = \gamma_{\Delta_c} + c - 1$.

Let the minimal generators of $\Delta_c$ be $\gamma_1 < \dots < \gamma_m$. We only need to prove $\gamma_m \le N'$, which guarantees that all generators lie in $\Gamma_{\leq N'}$.

Observe that the difference between any two distinct generators of $\Delta_c$ must belong to the set of gaps $\mathbb{N} \setminus \Gamma$. Since there are at most $c-1$ gaps, we have $\gamma_m - \gamma_1 \le c-1$. Because $\gamma_1 < \gamma_{\Delta_c}$, it follows that
\[
\gamma_m < \gamma_{\Delta_c} + c - 1 = N'.
\]
Thus $\gamma_1,\dots,\gamma_m \in \Gamma_{\leq N'}$, and the algorithm runs correctly with this truncation.
\end{proof}
\end{lemma}

By a similar argument, the same bound $\gamma_{\Delta_c} + c - 1$ is sufficient for Algorithm~\ref{algorithm tree}.

The situation for Algorithm~\ref{algorithm syz} is slightly different: the generators of $\operatorname{Syz}(\Delta)$ may exceed $\gamma_{\Delta_c} + c - 1$. Nevertheless, within Algorithm~\ref{algorithm motivic} we only need to know which generators of $\operatorname{Syz}(\Delta)$ are smaller than $\gamma_\Delta$, and $\gamma_\Delta \le \gamma_{\Delta_c} < \gamma_{\Delta_c} + c - 1$. Hence the bound obtained for Algorithm~\ref{algorithm min generators} also guarantees the correct execution of Algorithm~\ref{algorithm syz} inside the main algorithm.

\begin{theorem}\label{limite of effective bound}
Let $N$ be the minimal effective bound of $\Gamma$ required by Algorithm~\ref{algorithm motivic}. Then
\[
N \le (c-1)(\alpha_1 + 2).
\]
\begin{proof}
Let  $\Delta_{c}\in \mathscr{D}_{c}$ such that $\gamma_{\Delta_{c}}=\max \{\gamma_{\Delta}\}_{\Delta\in \mathscr{D}_{c}}$. Then $N\leq \gamma_{\Delta_{c}}+1+c$. Set $\Delta_{c-1}= \Delta_{c}\cup \{\gamma_{\Delta_{c}}\}$. If $\gamma_{\Delta_{c}}-\alpha_{1}\in \Gamma$, then $\gamma_{\Delta_{c}}-\alpha_{1}\notin \Delta_{c-1}$. Thus,  $\gamma_{\Delta_{c-1}}\geq \gamma_{\Delta_{c}}-\alpha_{1}$ implying $\gamma_{\Delta_{c}}-\gamma_{\Delta_{c-1}}\leq \alpha_{1}$. Continue this process, we find that there exists a $t$ such that $\gamma_{\Delta_{t}}-\alpha_{1}\notin \Gamma$, which implies  $\gamma_{\Delta_{t}}-\alpha_{1}\leq c-1$, i.e., $\gamma_{\Delta_{t}}\leq \alpha_{1}+c-1$. Since $t\geq 2$, we have $\gamma_{\Delta_{c}}\leq (c-t+1)\alpha_{1}+c-1\leq (c-1)\alpha_{1}+c-1$. Thus, $N\leq (c-1)(\alpha_{1}+2)$.\end{proof}
\end{theorem}

\section{Experimentation}\label{section Experimentation}
This section provides an example illustrating the computation of the motivic Hilbert zeta function for certain curve singularities. 

\begin{example}
Let $C$ be a curve singularity with local ring $\widehat{\mathcal{O}_{C,O}} = \mathbb{C}[[t^{4},t^{5},t^{6}]]$, whose valuation semigroup is monomial: $\Gamma = \Gamma_{4,3}$ in the sense of Proposition~\ref{monomial semigroup}. Explicitly, $\Gamma = \langle 4,5,6 \rangle$ with conductor $c = 8$. According to the bound given in Theorem~\ref{limite of effective bound}, we set $n = (c-1)(\alpha_{1}+2) = 7 \times 6 = 42$. This means we implement the truncated set $\Gamma_{\leq 42}$ in our computations using Algorithm~\ref{algorithm motivic}, where
\[
\Gamma_{\leq 42} = \{0, 4, 5, 6, 8\} \cup [9,42].
\]

We initialize with $\mathscr{D}_{1} = \{(4,5,6)\}$, where $(4,5,6)$ denotes the set $\{0, 4, 5, 6, 8\} \cup [9,42]$. Applying Algorithm~\ref{algorithm tree} yields the following collections of $\Gamma$-subsemimodules:

\noindent
\begin{align*}
\mathscr{D}_{2} &= \{(5, 6, 8), (4, 6), (4, 5)\},\\
\mathscr{D}_{3} &= \{(6, 8, 9), (5, 8), (5, 6), (4, 11)\},\\
\mathscr{D}_{4} &= \{(8, 9, 10, 11), (6, 9), (6, 8), (5, 12), (4)\},\\
\mathscr{D}_{5} &= \{(9, 10, 11, 12), (8, 10, 11), (8, 9, 11), (8, 9, 10), (6, 13), (5)\},\\
\mathscr{D}_{6} &= \{(10, 11, 12, 13), (9, 11, 12), (9, 10, 12), (9, 10, 11), (8, 11), (8, 10), (8, 9), (6)\},\\
\mathscr{D}_{7} &= \{(11, 12, 13, 14), (10, 12, 13), (10, 11, 13), (10, 11, 12), (9, 12), (9, 11), (9, 10), (8, 15)\},\\
\mathscr{D}_{8} &= \{(12, 13, 14, 15), (11, 13, 14), (11, 12, 14), (11, 12, 13), (10, 13), (10, 12), (10, 11), (9, 16), (8)\}.
\end{align*}

The motivic classes of the punctual Hilbert schemes in the Grothendieck ring are computed as follows:
\[
\begin{aligned}
[C^{[1]}] &= 1,\\
[C^{[2]}] &= 1 + \mathbb{L} + \mathbb{L}^{2},\\
[C^{[3]}] &= 1 + \mathbb{L} + 2\mathbb{L}^{2},\\
[C^{[4]}] &= 1 + \mathbb{L} + 2\mathbb{L}^{2} + \mathbb{L}^{3},\\
[C^{[5]}] &= 1 + \mathbb{L} + 2\mathbb{L}^{2} + 2\mathbb{L}^{3},\\
[C^{[6]}] &= 1 + \mathbb{L} + 2\mathbb{L}^{2} + 3\mathbb{L}^{3} + \mathbb{L}^{4},\\
[C^{[7]}] &= 1 + \mathbb{L} + 2\mathbb{L}^{2} + 3\mathbb{L}^{3} + \mathbb{L}^{4},\\
[C^{[8]}] &= 1 + \mathbb{L} + 2\mathbb{L}^{2} + 3\mathbb{L}^{3} + 2\mathbb{L}^{4}.
\end{aligned}
\]

From these we obtain the motivic Hilbert zeta function
\[
Z^{\mathrm{Hilb}}_{(C,O)}(q) = 1 + \sum_{\ell=1}^{7} [C^{[\ell]}] q^{\ell} + [C^{[8]}] (q^{8} + q^{9} + \cdots).
\]

Finally, take $\Delta = (8) \in \mathscr{D}_{8}$ and observe that $\gamma_{\Delta} = \max_{\Delta' \in \mathscr{D}_{c}} \gamma_{\Delta'} = 15$. We find
\[
N \le \gamma_{\Delta} + c - 1 = 22 \le 42 = (c-1)(\alpha_{1}+2),
\]
which is consistent with Theorem~\ref{limite of effective bound}.
\end{example}

\begin{example}
We consider $C=\{y^3=x^7\}$,  its semigroup is $\Gamma= \langle 3, 7 \rangle$ with conductor $c=12$. 

\[
\begin{aligned}
Zm_{(C,O)}^{\mathrm{Hilb}}(a^2, q^2) 
&= \sum_{\ell \geq 0} \sum_{m \geq 1} q^{2\ell} (1 - a^2)^{m-1} [C^{[\ell],m}] \\
&= \sum_{\ell \geq 0} \sum_{\Delta \in \mathscr{D}_\ell} \sum_{m \geq 1} q^{2\ell} (1 - a^2)^{m-1} [C^{[\Delta],m}].
\end{aligned}
\]

The computation of $Zm_{(C,O)}^{\mathrm{Hilb}}$ proceeds by evaluating the sum layer-wise with respect to the index $\ell$:\\

For $\ell=0$:
\[
\sum_{m \geq 1} q^{2\ell} (1 - a^2)^{m-1} [C^{[\ell],m}] = 1.
\]

For $\ell=1$:
\[
\sum_{m \geq 1} q^{2\ell} (1 - a^2)^{m-1} [C^{[\ell],m}] = q^2(1 - a^2).
\]

For $\ell=2$:
\[
\sum_{m \geq 1} q^{2\ell} (1 - a^2)^{m-1} [C^{[\ell],m}] = (1 + \mathbb{L}) q^4(1 - a^2).
\]

For $\ell=3$:
\[
\sum_{m \geq 1} q^{2\ell} (1 - a^2)^{m-1} [C^{[\ell],m}] = q^6(1 - a^2)^2 + \mathbb{L} q^6(1 - a^2) + \mathbb{L}^2 q^6.
\]

For $\ell=4$:
\[
\sum_{m \geq 1} q^{2\ell} (1 - a^2)^{m-1} [C^{[\ell],m}] = q^8(1 - a^2)^2 + (\mathbb{L} + 2\mathbb{L}^2) q^8(1 - a^2).
\]

For $\ell=5$:
\[
\sum_{m \geq 1} q^{2\ell} (1 - a^2)^{m-1} [C^{[\ell],m}] = (1 + \mathbb{L}) q^{10}(1 - a^2)^2 + (2\mathbb{L}^2 + \mathbb{L}^3) q^{10}(1 - a^2).
\]

For $\ell=6$:
\[
\sum_{m \geq 1} q^{2\ell} (1 - a^2)^{m-1} [C^{[\ell],m}] = (1 + \mathbb{L}) q^{12}(1 - a^2)^2 + (\mathbb{L}^2 + \mathbb{L}^3) q^{12}(1 - a^2) + \mathbb{L}^4 q^{12}.
\]

For $\ell=7$:
\[
\sum_{m \geq 1} q^{2\ell} (1 - a^2)^{m-1} [C^{[\ell],m}] = (1 + \mathbb{L} + 2\mathbb{L}^2) q^{14}(1 - a^2)^2 + (2\mathbb{L}^3 + \mathbb{L}^4) q^{14}(1 - a^2) + 2\mathbb{L}^4 q^{14}.
\]

For $\ell=8$:
\[
\sum_{m \geq 1} q^{2\ell} (1 - a^2)^{m-1} [C^{[\ell],m}] = (1 + \mathbb{L} + 2\mathbb{L}^2) q^{16}(1 - a^2)^2 + (2\mathbb{L}^3 + 3\mathbb{L}^4) q^{16}(1 - a^2).
\]

For $\ell=9$:
\[
\sum_{m \geq 1} q^{2\ell} (1 - a^2)^{m-1} [C^{[\ell],m}] = (1 + \mathbb{L} + 2\mathbb{L}^2 + \mathbb{L}^3) q^{18}(1 - a^2)^2 + (3\mathbb{L}^3 + 2\mathbb{L}^4) q^{18}(1 - a^2) + \mathbb{L}^5 q^{18}.
\]

For $\ell=10$:
\[
\sum_{m \geq 1} q^{2\ell} (1 - a^2)^{m-1} [C^{[\ell],m}] = (1 + \mathbb{L} + 2\mathbb{L}^2 + \mathbb{L}^3) q^{20}(1 - a^2)^2 + (3\mathbb{L}^3 + 2\mathbb{L}^4 + \mathbb{L}^5) q^{20}(1 - a^2) + \mathbb{L}^5 q^{20}.
\]

For $\ell=11$:
\[
\sum_{m \geq 1} q^{2\ell} (1 - a^2)^{m-1} [C^{[\ell],m}] = (1 + \mathbb{L} + 2\mathbb{L}^2 + \mathbb{L}^3) q^{22}(1 - a^2)^2 + (3\mathbb{L}^3 + 2\mathbb{L}^4 + 2\mathbb{L}^5) q^{22}(1 - a^2).
\]

For $\ell=12$:
\[
\sum_{m \geq 1} q^{2\ell} (1 - a^2)^{m-1} [C^{[\ell],m}] = (1 + \mathbb{L} + 2\mathbb{L}^2 + \mathbb{L}^3) q^{24}(1 - a^2)^2 + (3\mathbb{L}^3 + 2\mathbb{L}^4 + 2\mathbb{L}^5) q^{24}(1 - a^2) + \mathbb{L}^6 q^{24}.
\]

In summary, we have 
\begin{align*}
Zm_{(C,O)}^{\mathrm{Hilb}}(a^2, q^2) 
 &= \sum_{\ell=0}^{12}  \sum_{m \geq 1} q^{2\ell} (1 - a^2)^{m-1} [C^{[\ell],m}]   +  \sum_{\ell> 12}  \sum_{m \geq 1} q^{2\ell} (1 - a^2)^{m-1} [C^{[\ell],m}]  \\
& = 1 \\
& + (1-a^2) q^2 \\
& + (1+\mathbb{L})  (1-a^2) q^4\\
& + \left[ \mathbb{L}^2 + \mathbb{L}(1-a^2) + (1-a^2)^2 \right] q^6 \\
& + \left[ (\mathbb{L} + 2\mathbb{L}^2)(1-a^2) + (1-a^2)^2 \right] q^8 \\
& + \left[ (2\mathbb{L}^2 + \mathbb{L}^3)(1-a^2) + (1+\mathbb{L})(1-a^2)^2 \right] q^{10} \\
& + \left[ \mathbb{L}^4 + (\mathbb{L}^2 + \mathbb{L}^3)(1-a^2) + (1+\mathbb{L})(1-a^2)^2 \right] q^{12} \\
& + \left[ 2\mathbb{L}^4 + (2\mathbb{L}^3 + \mathbb{L}^4)(1-a^2) + (1+\mathbb{L} + 2\mathbb{L}^2)(1-a^2)^2 \right] q^{14} \\
& + \left[ (2\mathbb{L}^3 + 3\mathbb{L}^4)(1-a^2) + (1+\mathbb{L} + 2\mathbb{L}^2)(1-a^2)^2 \right] q^{16} \\
& + \left[ \mathbb{L}^5 + (3\mathbb{L}^3 + 2\mathbb{L}^4)(1-a^2) + (1+\mathbb{L} + 2\mathbb{L}^2 + \mathbb{L}^3)(1-a^2)^2 \right] q^{18} \\
& + \left[ \mathbb{L}^5 + (3\mathbb{L}^3 + 2\mathbb{L}^4 + \mathbb{L}^5)(1-a^2) + (1+\mathbb{L} + 2\mathbb{L}^2 + \mathbb{L}^3)(1-a^2)^2 \right] q^{20} \\
& + \left[ (3\mathbb{L}^3 + 2\mathbb{L}^4 + 2\mathbb{L}^5)(1-a^2) + (1+\mathbb{L} + 2\mathbb{L}^2 + \mathbb{L}^3)(1-a^2)^2 \right] q^{22} \\
& + \left[ \mathbb{L}^6 + (3\mathbb{L}^3 + 2\mathbb{L}^4 + 2\mathbb{L}^5)(1-a^2) + (1+\mathbb{L} + 2\mathbb{L}^2 + \mathbb{L}^3)(1-a^2)^2 \right]( q^{24}+q^{26}+\dots ). 
\end{align*}
\end{example}

\textbf{Acknowledgments.}
We especially grateful to Mounir Hajli, another supervisor of Wenhao Zhu, deeply for his instruction on mathematics and his sustained encouragement to the first author and long term useful discussions and suggestions in this topic. We thank Ilaria Rossinelli for several very helpful discussions. We thank André Belotto da Silva for useful conversations and suggestions. The first author being supported by the  funding from  the Guangzhou Elites Sponsorship Council under the Oversea Study Program of Guangzhou Elite Project. 


\printbibliography
\end{document}